\documentclass[10pt,reqno]{amsart}
\usepackage{amsmath,amssymb,bbm,mathrsfs,extpfeil,graphicx}
\usepackage[margin=1in]{geometry}
\usepackage{ifthen,latexsym,float,colortbl,caption,dsfont,arydshln,sidecap}
\usepackage[mathscr]{eucal}
\usepackage{rotating,stackengine,enumerate,xcolor}
\usepackage{hyperref}

\DeclareMathOperator*{\esssup}{ess\,sup}

\definecolor{black}{rgb}{0.0, 0.0, 0.0}
\definecolor{red}{rgb}{1.0, 0.5, 0.5}

\provideboolean{shownotes}
\setboolean{shownotes}{true}

\newcommand{\margnote}[1]{
	\ifthenelse{\boolean{shownotes}}
	{\marginpar{\raggedright\tiny\texttt{#1}}}
	{}
}
\newcommand{\hole}[1]{
	\ifthenelse{\boolean{shownotes}}
	{\begin{center} \fbox{ \rule {.25cm}{0cm} \rule[-.1cm]{0cm}{.4cm}
				\parbox{.85\textwidth}{\begin{center} \texttt{#1}\end{center}} \rule
				{.25cm}{0cm}}\end{center}} {} }

\graphicspath{{pics/}}

\title[Global existence and semi-classical limit of Fermi--Dirac Fokker--Planck equation]
{Global existence and semi-classical limit for the spatially inhomogeneous Fermi--Dirac Fokker--Planck equation}

\author[Hyun]{Ju-Hwan Hyun}
\address[Ju-Hwan Hyun]{\newline Department of Mathematics,
\newline Yonsei University, 50 Yonsei-Ro, Seodaemun-Gu, Seoul 03722, Republic of Korea}
\email{jhhyun@yonsei.ac.kr}

\numberwithin{equation}{section}

\newtheorem{theorem}{Theorem}[section]
\newtheorem{lemma}{Lemma}[section]

\newtheorem{definition}{Definition}[section]

\makeatletter
\def\@subsubsection#1{\par
  \addpenalty\@secpenalty
  \vskip 1.5ex plus 1ex minus .2ex
  \refstepcounter{subsubsection}
  \addcontentsline{toc}{subsubsection}{\protect\numberline{\thesubsubsection}#1}
  \noindent{\normalfont\bfseries\thesubsubsection\quad#1\par\nobreak}
  \vskip 0.5ex \@afterheading}
\makeatother

\newcommand{\lt}{\left}
\newcommand{\rt}{\right}

\newcommand{\mN}{\mathbb{N}}

\newcommand{\mR}{\mathbb{R}}

\newcommand{\mcC}{\mathcal{C}}
\newcommand{\mcD}{\mathcal{D}}

\newcommand{\mcF}{\mathcal{F}}

\newcommand{\mcH}{\mathcal{H}}

\newcommand{\msF}{\mathscr{F}}

\newcommand{\fh}{f_\hbar}

\newcommand{\inttr}{\iint_{\mR^d \times \mR^d}}

\newcommand{\rd}{\textnormal{d}}

\newcommand{\rdx}{\textnormal{d}x}

\newcommand{\rdt}{\textnormal{d}t}

\newcommand{\rdy}{\textnormal{d}y}

\newcommand{\rds}{\textnormal{d}s}
\newcommand{\rdp}{\textnormal{d}p}
\newcommand{\rdq}{\textnormal{d}q}

\newcommand{\supp}{{\rm supp}}

\begin{document}

\allowdisplaybreaks
\date{\today}

\begin{abstract}
    In this paper, we study the global existence and the semi-classical limit of the spatially inhomogeneous Fermi--Dirac Fokker--Planck equation.
    We first establish the existence and uniqueness of a global weak solution for initial data with finite mass and second velocity moment satisfying the Pauli bound.
    When the initial datum also has finite spatial second moment and finite free energy, we prove that the solution satisfies the Fermi--Dirac free energy inequality.
    Next, as the quantum parameter $\hbar$ goes to zero, we prove that the solution of the Fermi--Dirac Fokker--Planck equation converges strongly to the classical Fokker--Planck solution in $\mcC([0,T];L^1_j(\mR^d\times\mR^d))$ for all $j\in[0,2)$.
    We also obtain a quantitative convergence rate in $\hbar$ when both equations share the same bounded initial datum.
\end{abstract}

\maketitle

\section{Introduction}

The kinetic description of quantum gases was first proposed by Nordheim \cite{Nor28}, and Uehling and Uhlenbeck studied transport phenomena in Fermi--Dirac and Bose--Einstein statistics \cite{UU33}.
Later, Kaniadakis and Quarati introduced a nonlinear kinetic Fokker--Planck model with quantum statistical effects \cite{KQ93, KQ94}.
Such quantum effects are included through statistical factors in the collision terms of quantum kinetic equations.
For fermions, the Pauli exclusion principle gives an upper bound on the particle density in the phase space.

We consider the spatially inhomogeneous Fermi--Dirac Fokker--Planck equation for the particle density $\fh=\fh(x,p,t)\ge 0$ with $(x,p,t)\in\mR^d\times\mR^d\times[0,\infty)$:
\begin{align}\label{eq:quantum}
    \partial_t \fh+p\cdot \nabla_x\fh=\nabla_p\cdot\big(\nabla_p \fh+p\fh(1-\hbar \fh)\big),
\end{align}
with initial data $\fh(x,p,0)=f_{\hbar,0}(x,p)$.
Here $\hbar>0$ is the quantum parameter and the factor $1-\hbar \fh$ comes from the Pauli exclusion principle, which gives the bound $0\le \fh\le 1/\hbar$.
As discussed in \cite{CHHglobal}, one can formally derive \eqref{eq:quantum} from the quantum Landau equation by the grazing collision limit \cite{Vil98} under radial symmetry and normalization assumptions.

Several results have been obtained on the existence of solutions for quantum kinetic models.
For the spatially homogeneous case, quantum Boltzmann equations were studied in \cite{Dol94, EMV05, Lu01, LW03}, and the existence and large-time behavior for the Landau--Fermi--Dirac equation were studied in \cite{ABDL22, ABL21, BL04, GGZ22}.
For the spatially homogeneous Fermi--Dirac Fokker--Planck equation, Carrillo, Lauren\c{c}ot, and Rosado obtained global weak solutions and studied their large-time behavior by using the entropy and its dissipation \cite{CLR09}, while the Bose--Einstein case was treated in \cite{CCLR16, Tos12}.
In the spatially inhomogeneous case, global solutions near steady states were studied for the quantum Fokker--Planck equation in \cite{CHHglobal, LZ15, NS07} and for the quantum BGK model in \cite{BY20}, while global solutions for the Landau--Fermi--Dirac equation were obtained in \cite{Sam24}.

Recently, Toshpulatov \cite{Tos26} studied a model in which the spatial transport term takes the nonlinear form $p\cdot\nabla_x(f-f^2)$, and obtained global $BV$ entropy solutions by adding a spatial viscosity term $\epsilon\Delta_x f_\epsilon$ under regularity and moment assumptions on the initial data.
In our equation \eqref{eq:quantum}, the spatial transport term $p\cdot\nabla_x\fh$ is linear, so we directly use the Green's function of the linear kinetic Fokker--Planck equation \cite{VO90} without spatial viscosity.
From the estimates on the Green's function, we first construct global solutions for regular initial data by using the Pauli bound $0\le\fh\le 1/\hbar$, conservation of mass, and velocity moment bounds.
Also, since the spatial transport term is linear, we obtain the $L^1$-contraction property without spatial derivative bounds on the solutions, and we approximate the initial data in $L^1_2(\mR^d\times\mR^d)$ to construct a unique global weak solution.

We now introduce the weighted Lebesgue spaces and the definition of weak solutions used throughout this work.
For $k\in[1,\infty)$ and $m\ge0$, let $L^k_m(\mR^d\times\mR^d)$ be the weighted Lebesgue space with the norm
$$\|g\|_{L^k_m}:=\lt(\inttr(1+|p|^m)^k|g(x,p)|^k\,\rdx\rdp\rt)^{\frac{1}{k}},$$
and for $m=0$, we write $L^k_{x,p}$ and $\|g\|_{L^k_{x,p}}$.

\begin{definition}[Weak solution]\label{def:weak solution}
    Let $T>0$ and $f_{\hbar,0}\in L^1_2(\mR^d\times\mR^d)$ satisfy $0\le f_{\hbar,0}\le 1/\hbar$.
    A function $\fh$ is called a weak solution of \eqref{eq:quantum} on $[0,T]$ with initial data $f_{\hbar,0}$ if it satisfies
    $$\fh\in \mcC([0,T];L^1_1(\mR^d\times\mR^d))\cap L^\infty(0,T;L^1_2(\mR^d\times\mR^d)\cap L^\infty(\mR^d\times\mR^d)).$$
    Moreover, $\fh$ satisfies $0\le\fh\le1/\hbar$ and $\nabla_p\fh\in L^2(0,T;L^2(\mR^d\times\mR^d))$, and the following weak formulation holds for every $\phi\in\mcC_c^\infty([0,T)\times\mathbb R^d\times\mathbb R^d)$:
    \begin{align}\label{eq:weak form 1}
    \begin{split}
        &-\int_0^T\inttr \fh(\partial_t\phi+p\cdot\nabla_x\phi+\Delta_p\phi-p\cdot\nabla_p\phi)\,\rdx\rdp\rdt\cr
        &\quad=\inttr f_{\hbar,0}\phi(x,p,0)\,\rdx\rdp+\int_0^T\inttr \hbar p\fh^2\cdot\nabla_p\phi\,\rdx\rdp\rdt.
    \end{split}
    \end{align}
\end{definition}

Also, for $0\le g\le 1/\hbar$, we put $0\log0=0$ and define the Fermi--Dirac free energy $\mcH_\hbar$ and the dissipation $\mcD_\hbar$ as
\begin{align*}
    \mcH_\hbar(g)&:=\inttr\lt(g\log g+\frac{1}{\hbar}(1-\hbar g)\log(1-\hbar g)+\frac{|p|^2}{2}g\rt)\,\rdx\rdp,\cr
    \mcD_\hbar(g)&:=\inttr g(1-\hbar g)\lt|\nabla_p\lt(\log\frac{g}{1-\hbar g}+\frac{|p|^2}{2}\rt)\rt|^2\,\rdx\rdp,
\end{align*}
where the integrand in $\mcD_\hbar(g)$ is zero at $(x,p)$ when $g(x,p)=0$ or $g(x,p)=1/\hbar$.

\begin{theorem}[Global weak solution and free energy inequality]\label{theorem:exist}
    Assume that the initial datum $f_{\hbar,0}\in L^1_2(\mR^d\times\mR^d)$ of \eqref{eq:quantum} satisfies $0\le f_{\hbar,0}\le 1/\hbar$.
    Then, for any $T>0$, there is a unique global weak solution $\fh$ of \eqref{eq:quantum} on $[0,T]$ in the sense of Definition \ref{def:weak solution}, and for all $0\le t\le T$,
    $$0 \le \fh(x,p,t) \le \frac{1}{\hbar}\quad\text{and}\quad \|\fh(t)\|_{L^1_{x,p}}=\|f_{\hbar,0}\|_{L^1_{x,p}}.$$
    In addition, if $f_{\hbar,0}$ satisfies
    $$\inttr |x|^2f_{\hbar,0}(x,p)\,\rdx\rdp<\infty\quad\text{and}\quad\mcH_\hbar(f_{\hbar,0})<\infty,$$
    then for all $0\le t\le T$, we have $\inttr |x|^2\fh(x,p,t)\,\rdx\rdp<\infty$ and
    \begin{align}\label{ineq:main_free_energy}
        \mcH_\hbar(\fh(t))+\int_0^t\mcD_\hbar(\fh(s))\,\rds\le \mcH_\hbar(f_{\hbar,0}).
    \end{align}
\end{theorem}

Our second main result is the semi-classical limit as $\hbar\to0$.
By Bohr's correspondence principle, quantum models are expected to recover classical models when the quantum parameter $\hbar$ goes to zero.
If $\hbar$ tends to zero, the nonlinear drift term $\nabla_p\cdot(p\fh(1-\hbar\fh))$ formally converges to the linear drift term $\nabla_p\cdot(pf)$, and \eqref{eq:quantum} is expected to converge to the classical linear Fokker--Planck equation:
\begin{align}\label{eq:classic}
    \partial_t f+p\cdot\nabla_xf=\nabla_p\cdot(\nabla_pf+pf),
\end{align}
with initial data $f(x,p,0)=f_0(x,p)$.

To pass to the limit $\hbar\to0$, we need to control the nonlinear term $\hbar p\fh^2$.
From the Pauli bound $0\le\fh\le1/\hbar$, we only have $\hbar\fh^2\le\fh$, which does not vanish as $\hbar\to0$.
Also, by the Pauli bound and the $L^r$ estimate in Section \ref{sec:existence} (letting $r\to\infty$ in \eqref{ineq:Lr_bound}), we obtain the following $L^\infty$ estimate:
$$\|f_\hbar(t)\|_{L^\infty_{x,p}}\le\min\lt\{\frac1\hbar,e^{dt}\|f_{\hbar,0}\|_{L^\infty_{x,p}}\rt\}.$$
This bound is still not uniform as $\hbar\to0$, since $\|f_{\hbar,0}\|_{L^\infty_{x,p}}$ may depend on $\hbar$.
Thus, the $L^\infty$ and $L^2$ bounds of $\fh$ used for fixed $\hbar>0$ cannot be directly applied to pass to the limit.
Instead, we approximate the initial data by a fixed regular function to avoid the $\hbar$-dependence, and then combine the $L^1$-contraction property with the uniform velocity second moment bound.

In recent years, the semi-classical limit as $\hbar\to0$ has been studied for several other quantum kinetic equations.
For the spatially homogeneous quantum Boltzmann equation, He, Lu, and Pulvirenti \cite{HLP21} showed weak convergence to the classical Fokker--Planck--Landau equation, and He, Lu, Pulvirenti, and Zhou \cite{HLPZ24} later obtained a local-in-time asymptotic expansion in weighted Sobolev spaces when both equations have the same initial datum
In the spatially inhomogeneous setting, Sampaio \cite{Sam26} proved convergence up to a subsequence for the Landau--Fermi--Dirac equation with Coulomb potential to a renormalized solution of the classical Landau equation with a defect measure.
For the spatially inhomogeneous Fermi--Dirac BGK model, Choi, Hwang, and Song \cite{CHS26} showed strong convergence up to a subsequence under uniform moment and entropy bounds.

\begin{theorem}[Semi-classical limit and convergence rate]\label{theorem:semi}
    Let $f_0\in L^1_2(\mR^d\times\mR^d)$ be a nonnegative initial datum for \eqref{eq:classic}, and let $f_{\hbar,0}\in L^1_2(\mR^d\times\mR^d)$ be the initial data for \eqref{eq:quantum}.
    Suppose that
    \begin{align}\label{assu:hbar initial}
        0\le f_{\hbar,0}(x,p)\le\frac{1}{\hbar},\quad \lim_{\hbar\to0}\|f_{\hbar,0}-f_0\|_{L^1_{x,p}}=0,
    \end{align}
    and
    \begin{align}\label{assu:moment}
        \sup_{0<\hbar\le1}\inttr(1+|p|^2)f_{\hbar,0}(x,p)\,\rdx\rdp<\infty.
    \end{align}
    Let $f_\hbar$ and $f$ be the unique global weak solutions of \eqref{eq:quantum} and \eqref{eq:classic}, respectively.
    Then, for any $j\in[0,2)$ and $T>0$, we have
    $$f_\hbar\to f\quad\text{strongly in}\quad\mcC([0,T];L^1_j(\mR^d\times\mR^d))\quad\text{as }\hbar\to0.$$
    Also, if $f_0\in L^\infty(\mR^d\times\mR^d)$, then for $g_\hbar$ solving \eqref{eq:quantum} with $g_\hbar(x,p,0)=f_0(x,p)$ and sufficiently small $\hbar>0$, there is a constant $C>0$ depending only on $d$, $T$, $\|f_0\|_{L^1_2}$, and $\|f_0\|_{L^\infty_{x,p}}$ such that for all $j\in[0,2)$,
    \begin{align}\label{ineq:quant_rate}
        \sup_{0\le t\le T}\|g_\hbar(t)-f(t)\|_{L^1_j}\le C\hbar^{\min\{1,\,2-j\}}.
    \end{align}
\end{theorem}

Theorem \ref{theorem:semi} justifies the semi-classical limit for the spatially inhomogeneous Fermi--Dirac Fokker--Planck equation \eqref{eq:quantum} with linear spatial transport.
Unlike the previous spatially inhomogeneous results on the Landau--Fermi--Dirac and BGK models in \cite{CHS26, Sam26}, where convergence was obtained up to a subsequence, we prove strong convergence in $\mcC([0,T];L^1_j(\mR^d\times\mR^d))$ for the whole family $\fh$ as $\hbar\to0$ without extracting a subsequence.
In addition, when $f_0\in L^\infty(\mR^d\times\mR^d)$ is taken as the common initial datum for both equations, we combine the $L^\infty$ and velocity moment bounds to derive the explicit $O(\hbar^{\min\{1,\,2-j\}})$ convergence rate \eqref{ineq:quant_rate} in $L^1_j(\mR^d\times\mR^d)$ for all $j\in[0,2)$.

The rest of this paper is organized as follows.
In Section \ref{sec:existence}, we review the Green's function of the linear kinetic Fokker--Planck equation and the global existence of solutions for regular initial data, and then we prove Theorem \ref{theorem:exist}.
In Section \ref{sec:semi}, we state the basic properties of the classical solution and prove Theorem \ref{theorem:semi}.

\section{Global weak solutions and free energy inequality}\label{sec:existence}

\subsection{Properties of the Green's function and local existence}

Let $\mcF(x,y,p,q,t)$ be the two-point Green's function of the linear Fokker--Planck equation \eqref{eq:classic}.
The explicit formula of $\mcF$ is given in \cite{VO90} as follows:
\begin{align}\label{def:Green}
\begin{split}
    \mcF(x,y,p,q,t) &:= \lt[\frac{e^{t}}{4\pi (\tilde{a}(t)t-(a(t))^2)^{1/2}}\rt]^{d} \cr
    &\quad\times\exp{\lt\{-\frac{|x-y+p-q|^2}{4t}-\frac{|a(t)(x-y+p-q)+t(q-pe^t)|^2}{4t(\tilde{a}(t)t-(a(t))^2)}\rt\}},
\end{split}
\end{align}
where $a(t)=e^{t}-1$ and $\tilde{a}(t)=\frac{e^{2t}-1}{2}.$
As shown in \cite{VO90}, for $t-s>0$, $\mcF$ satisfies
\begin{align}\label{prop}
    \partial_t \mcF(x,y,p,q,t-s)+p\cdot\nabla_x\mcF(x,y,p,q,t-s)=\nabla_p\cdot(\nabla_p\mcF(x,y,p,q,t-s)+p\mcF(x,y,p,q,t-s)),
\end{align}
and for $t>0$,
\begin{align}\label{prop:Green_mass}
    \inttr \mcF(x,y,p,q,t)\,\rdy\rdq=e^{dt}\quad\text{and}\quad \inttr \mcF(x,y,p,q,t)\,\rdx\rdp=1.
\end{align}

Now, we define the linear operators $\msF[g](x,p,t)$ and $\msF^q[h](x,p,t)$, acting on functions $g:\mR^d\times\mR^d\to\mR$ and $h:\mR^d\times\mR^d\to\mR^d$ as
\begin{align*}
    \msF[g](x,p,t)&:=\inttr \mcF(x,y,p,q,t)g(y,q)\,\rdy\rdq, \cr
    \msF^q[h](x,p,t)&:=\inttr \nabla_q\mcF(x,y,p,q,t)\cdot h(y,q)\,\rdy\rdq.
\end{align*}
The estimates below follow from direct Gaussian kernel calculations and Young's inequality by the same argument as in \cite[Lemmas 4.1--4.3]{Tos26} (see also \cite{CLR09}), and we state them here for our Green's function $\mcF$.

\begin{lemma}\label{Lem:F bound}
    Let $1\le l\le k\le\infty$, $m=0\ \text{or}\ m\ge1$, and $T>0$.
    Then there are constants $C_1 = C_1(k,m,d,T) > 0$ and $C_2 = C_2(k,l,m,d,T) > 0$ such that
    \begin{align*}
        \|\msF[g](t)\|_{L^k_m}\le C_1 e^{dt}\|g\|_{L^k_m}, \quad \text{for all}\ t\in (0,T],
    \end{align*}
    and
    \begin{align}\label{ineq:F 2}
        \|\msF^{q}[h](t)\|_{L^k_m}\le \frac{C_2 e^{d\lt(1+\frac{1}{k}-\frac{1}{l}-\frac{1}{d}\rt)t}}{t^{2d\lt(\frac{1}{l}-\frac{1}{k}\rt)+\frac{1}{2}}}\|h\|_{L^l_m}, \quad \text{for all}\ t\in (0,T].
    \end{align}
    Also, the same bound \eqref{ineq:F 2} holds when $\nabla_q\mcF$ is replaced by $\nabla_p\mcF$.
\end{lemma}

To see how the time singularity in \eqref{ineq:F 2} arises, we examine the explicit formula \eqref{def:Green}.
When we estimate $\msF^q[h]$, differentiating the two terms in the exponent of $\mcF$ with respect to $q$ (or $p$) brings down factors of order $t^{-1/2}$ and 
$$\frac{|a(t)-t|}{\sqrt{t(\tilde{a}(t)t-(a(t))^2)}}\quad \lt(\text{or}\ \frac{|a(t)-te^t|}{\sqrt{t(\tilde{a}(t)t-(a(t))^2)}}\rt)$$
from the Gaussian kernel, while the prefactor in \eqref{def:Green} contributes $\big(\tilde{a}(t)t-(a(t))^2\big)^{-\frac{d}{2}(\frac{1}{l}-\frac{1}{k})}$ via Young's inequality.
By Taylor's theorem as $t\to 0^+$, the $O(t^2)$ and $O(t^3)$ terms in $\tilde{a}(t)t$ and $(a(t))^2$ cancel out:
$$a(t)-t = \frac{1}{2}t^2 + O(t^3), \quad a(t)-te^t = -\frac{1}{2}t^2 + O(t^3), \quad \tilde{a}(t)t-(a(t))^2 = \frac{1}{12}t^4 + O(t^5).$$
Hence, the derivative of the exponent gives $t^{-1/2}$ and $\big(\tilde{a}(t)t-(a(t))^2\big)^{-1/2}$ is of order $t^{-2}$ as $t\to 0^+$, so we get the time singularity $t^{-2d(\frac{1}{l}-\frac{1}{k})-\frac{1}{2}}$ in \eqref{ineq:F 2}.
On the other hand, \cite{Tos26} adds a spatial viscosity term $\epsilon\Delta_x f_\epsilon$ to treat the nonlinear transport $p\cdot\nabla_x(f-f^2)$.
Because of $\epsilon\Delta_x$, the determinant factor inside the square root of the Green's function in \cite{Tos26} has leading order $\epsilon t^2$ rather than $t^4$ near $t=0$.
For fixed $\epsilon>0$, this changes the prefactor order from $t^{-2d}$ to $t^{-d}$ and gives the weaker singularity $t^{-d(\frac{1}{l}-\frac{1}{k})-\frac{1}{2}}$ in \cite[Lemmas 4.2 and 4.3]{Tos26}.

Using the Green's function $\mcF$, we write \eqref{eq:quantum} in the integral form
\begin{align}\label{eq:int}
\begin{split}
     \fh(x,p,t) &= \msF[f_{\hbar,0}](x,p,t)+\int^t_0\msF^q\big[\hbar q(\fh(s))^2\big](x,p,t-s)\,\rds\cr
     &=\inttr \mcF(x,y,p,q,t)f_{\hbar,0}(y,q)\,\rdy\rdq \cr
     &\quad+\int^t_0\inttr\nabla_q\mcF(x,y,p,q,t-s)\cdot\big[\hbar q(\fh(y,q,s))^2\big]\,\rdy\rdq\rds.
\end{split}
\end{align}
To apply a fixed-point argument to \eqref{eq:int}, we set up the function spaces for the local solution.
For $k > 4d$ and $m \ge 1$, let $\Upsilon$ be the space
$$\Upsilon := L^\infty(\mR^d\times\mR^d) \cap L^1_1(\mR^d\times\mR^d) \cap L^k_m(\mR^d\times\mR^d),$$
equipped with the norm 
$$\|g\|_{\Upsilon} := \max\lt\{ \|g\|_{L^\infty_{x,p}}, \|g\|_{L^1_1}, \|g\|_{L^k_m} \rt\}.$$
For any $T>0$, we define $\Upsilon_T$ as
$$\Upsilon_T := \mcC([0,T]; L^1_1(\mR^d\times\mR^d) \cap L^k_m(\mR^d\times\mR^d)) \cap L^\infty(0,T; L^\infty(\mR^d\times\mR^d)),$$
with the norm
$$\|g\|_{\Upsilon_T} := \max\lt\{ \esssup_{0 < t \le T} \|g(t)\|_{L^\infty_{x,p}}, \sup_{0 \le t \le T} \|g(t)\|_{L^1_1}, \sup_{0 \le t \le T} \|g(t)\|_{L^k_m} \rt\}.$$
Also, to get the regularity of the approximated solution with respect to $p$, we consider the space
$$\Lambda_T:=\big\{g\in\Upsilon_T: \nabla_p g(t)\in L^1_1(\mR^d\times\mR^d)\cap L^k_m(\mR^d\times\mR^d)\ \text{for}\ 0<t\le T \ \text{and}\ \|g\|_{\Lambda_T}<\infty\big\}$$
with the weighted norm
$$\|g\|_{\Lambda_T} := \max\lt\{ \|g\|_{\Upsilon_T}, \esssup_{0 < t \le T} \sqrt{t} \|\nabla_p g(t)\|_{L^1_1}, \esssup_{0 < t \le T} \sqrt{t} \|\nabla_p g(t)\|_{L^k_m} \rt\}.$$
For $\fh \in \Lambda_T$, its derivative with respect to $p$ is given by
\begin{align*}
    \nabla_p\fh(x,p,t) &= \inttr \nabla_p\mcF(x,y,p,q,t)f_{\hbar,0}(y,q)\,\rdy\rdq \cr
    &\quad - \int^t_0\inttr \nabla_p\mcF(x,y,p,q,t-s)\nabla_q\cdot\big[\hbar q(\fh(y,q,s))^2\big]\,\rdy\rdq\rds.
\end{align*}
To estimate the nonlinear terms in \eqref{eq:int} and $\nabla_p\fh$ by Lemma \ref{Lem:F bound} with $l = k/2$ for the $L^k_m$ estimates and $l = k$ for the $L^\infty$ estimate, the time integral is finite when
$$\frac{2d}{k}+\frac{1}{2} < 1, \quad\text{that is,}\quad k > 4d.$$
Then, by using Lemma \ref{Lem:F bound} and the fixed-point argument in \cite{CLR09, Tos26}, we obtain the following local existence result, and Fubini's theorem with \eqref{prop} gives \eqref{eq:weak form 1}.

\begin{theorem}\label{thm:local}
    Let $k>4d$, $m \ge 1$, and $f_{\hbar,0} \in \Upsilon$.
    Then there exists $T > 0$ depending on $\hbar$ and $\|f_{\hbar,0}\|_{\Upsilon}$ such that \eqref{eq:quantum} has a unique solution $\fh \in \Lambda_T$ on $[0,T]$ satisfying the integral equation \eqref{eq:int}.
\end{theorem}

\subsection{Global existence and uniqueness}\label{sec:global}

In this subsection, we prove the global existence and uniqueness part of Theorem \ref{theorem:exist}.
Let us define a non-increasing function $\zeta \in \mcC^\infty_c([0,\infty))$ such that $\zeta(z)=1$ for $0\le z\le 1$ and $\zeta(z)=0$ for $z\ge 2$ with $|\zeta'|,|\zeta''|<2$, and define the cut-off functions $\zeta_N(x), \zeta_R(p)\in\mcC^\infty_c(\mR^d)$ as
$$\zeta_N(x)=\zeta\lt(\frac{|x|}{N}\rt),\quad \zeta_R(p)=\zeta\lt(\frac{|p|}{R}\rt)\quad \text{for $x,p\in\mR^d$ and $N,R\ge1.$}$$
Then we have
$$|\nabla_x\zeta_N(x)|\le \frac{C}{N}, \quad |\Delta_x\zeta_N(x)|\le \frac{C}{N^2}, \quad |\nabla_p\zeta_R(p)|\le \frac{C}{R}, \quad |\Delta_p\zeta_R(p)|\le \frac{C}{R^2}.$$

For initial data satisfying $0\le f_{\hbar,0}\le 1/\hbar$, the corresponding local solution $\fh\in\Lambda_T$ preserves the Pauli bound, conserves mass, and satisfies the velocity moment bounds by using the cut-off functions $\zeta_N(x)$ and $\zeta_R(p)$ in the same way as in \cite[Lemma 3.1 and Sections 4.4--4.5]{Tos26} (see also \cite{CLR09}).
In addition, for any $1 < r < \infty$, multiplying \eqref{eq:quantum} by $\fh^{r-1}\zeta_N(x)\zeta_R(p)$, integrating by parts, and letting $N,R\to\infty$ yield
\begin{align}\label{eq:Lr_identity}
\begin{split}
    \frac{1}{r}\frac{\rd}{\rdt}\inttr \fh^r\,\rdx\rdp + \frac{4(r-1)}{r^2}\inttr |\nabla_p \fh^{r/2}|^2\,\rdx\rdp &= \frac{d(r-1)}{r}\inttr \fh^r\,\rdx\rdp \cr
    &\quad - \frac{d(r-1)\hbar}{r+1}\inttr \fh^{r+1}\,\rdx\rdp.
\end{split}
\end{align}
From $-\frac{d(r-1)\hbar}{r+1}\inttr \fh^{r+1}\,\rdx\rdp \le 0$ in \eqref{eq:Lr_identity}, Gr\"onwall's inequality gives the $L^r$-estimate with the factor $e^{\frac{d(r-1)}{r}t}$ independent of $\hbar$, which is uniform in $\hbar$ on $[0,T]$ when the initial datum is fixed.
Also, by the convex approximation argument in \cite{CLR09, Tos26} and the linearity of $p\cdot\nabla_x\fh$, we get the $L^1$-contraction property without $x$-derivative bounds on $\fh$.

\begin{lemma}\label{Lem:reg_global}
    Suppose $m \ge 1$ and $k > 4d$. For $f_{\hbar,0} \in L^1_{mk}(\mR^d \times \mR^d)$ with $0 \le f_{\hbar,0} \le 1/\hbar$ and $T > 0$, \eqref{eq:quantum} has a unique global solution $\fh \in \Lambda_T$, and for $0\le t\le T$:
    \begin{enumerate}
        \item[(i)] $0 \le \fh(x,p,t) \le \frac{1}{\hbar}$ and $\|\fh(t)\|_{L^1_{x,p}} = \|f_{\hbar,0}\|_{L^1_{x,p}}.$
        \item[(ii)] All $2\gamma$-moments of $\fh(t)$ for $0 \le \gamma \le mk/2$ are bounded by a polynomial of $t$ depending only on $d$ and the initial moments. In particular, when $\gamma=1$, we have
        \begin{align}\label{ineq:second_moment}
            \inttr |p|^2\fh(x,p,t)\,\rdx\rdp\le \inttr |p|^2f_{\hbar,0}(x,p)\,\rdx\rdp+2d\|f_{\hbar,0}\|_{L^1_{x,p}}t.
        \end{align}
        \item[(iii)] For any $1 < r < \infty$, we have the $L^r$-estimate:
        \begin{align}\label{ineq:Lr_bound}
            \|\fh(t)\|_{L^r_{x,p}} \le e^{\frac{d(r-1)}{r}t}\|f_{\hbar,0}\|_{L^r_{x,p}}.
        \end{align}
        When $r=2$, we also have the $L^2$-estimate:
        \begin{align}\label{ineq:L2 energy}
            \frac{1}{2}\|\fh(t)\|_{L^2_{x,p}}^2 + \int_0^t \|\nabla_p \fh(s)\|_{L^2_{x,p}}^2\,\rds \le \frac{1}{2}e^{dt}\|f_{\hbar,0}\|_{L^2_{x,p}}^2.
        \end{align}
        \item[(iv)] More generally, if $\fh$ and $g_\hbar$ are weak solutions to \eqref{eq:quantum} on $[0,T]$ in the sense of Definition \ref{def:weak solution} with initial data $f_{\hbar,0}$ and $g_{\hbar,0}$, respectively, then
        \begin{align}\label{ineq:contraction}
            \|\fh(t) - g_\hbar(t)\|_{L^1_{x,p}} \le \|f_{\hbar,0} - g_{\hbar,0}\|_{L^1_{x,p}}.
        \end{align}
    \end{enumerate}
\end{lemma}

From (i) and (ii) of Lemma \ref{Lem:reg_global}, we see that $\|\fh(t)\|_{\Upsilon}$ is bounded on $[0,T]$ for any $T>0$, so the local solution in Theorem \ref{thm:local} extends globally in time.
We now prove the global existence and uniqueness part of Theorem \ref{theorem:exist} by using \eqref{ineq:contraction}.

\begin{proof}[Proof of Theorem \ref{theorem:exist} (Global existence and uniqueness)]
    Fix $m \ge 1$ and $k > 4d$. 
    For the initial data $f_{\hbar,0} \in L^1_2(\mR^d\times\mR^d)$ with $0\le f_{\hbar,0}\le \frac{1}{\hbar}$, we define the truncated initial data for $n\in\mN$ as
    $$f_{\hbar,0,n}(x,p) := f_{\hbar,0}(x,p)\mathbbm{1}_{\{|x|\le n, |p|\le n\}}.$$
    Note that since $f_{\hbar,0,n}(x,p)=0$ for $|p|>n$, we have
    $$\inttr (1+|p|^{mk})f_{\hbar,0,n}(x,p)\,\rdx\rdp \le (1+n^{mk})\|f_{\hbar,0}\|_{L^1_{x,p}} < \infty,$$
    and thus $f_{\hbar,0,n} \in L^1_{mk}(\mR^d\times\mR^d)$.
    Also, since $(1+|p|^2)|f_{\hbar,0,n} - f_{\hbar,0}| \le (1+|p|^2)f_{\hbar,0}$, we have
    \begin{align}\label{conv:init n}
        \lim_{n\to\infty}\|f_{\hbar,0,n} - f_{\hbar,0}\|_{L^1_2} = 0.
    \end{align}

    From Lemma \ref{Lem:reg_global} and \eqref{ineq:L2 energy}, we obtain the global solution $f_{\hbar,n} \in \Lambda_T$ of \eqref{eq:quantum} for $f_{\hbar,0,n}$ with the bounds:
    \begin{align}\label{bound:unif fn 1}
        0 \le f_{\hbar,n}(x,p,t) \le \frac{1}{\hbar}, \quad \|f_{\hbar,n}(t)\|_{L^1_{x,p}} = \|f_{\hbar,0,n}\|_{L^1_{x,p}} \le \|f_{\hbar,0}\|_{L^1_{x,p}},
    \end{align}
    and
    \begin{align}\label{bound:unif fn 2}
    \begin{split}
        \sup_{0\le t\le T}\inttr (1+|p|^2) f_{\hbar,n}(x,p,t)\,\rdx\rdp \le \inttr &(1+|p|^2) f_{\hbar,0}(x,p)\,\rdx\rdp + 2d\|f_{\hbar,0}\|_{L^1_{x,p}}T \le C_T, \cr
        \sup_{0\le t\le T}\frac{1}{2}\|f_{\hbar,n}(t)\|_{L^2_{x,p}}^2 +\int_0^T\|\nabla_p &f_{\hbar,n}(s)\|_{L^2_{x,p}}^2\,\rds \le e^{dT}\|f_{\hbar,0}\|_{L^2_{x,p}}^2\le C_{T,\hbar}.
    \end{split}
    \end{align}
    By \eqref{ineq:contraction}, for any $n, j \in \mN$, we have
    $$\sup_{0\le t\le T}\|f_{\hbar,n}(t) - f_{\hbar,j}(t)\|_{L^1_{x,p}} \le \|f_{\hbar,0,n} - f_{\hbar,0,j}\|_{L^1_{x,p}}.$$
    To pass to the limit in the nonlinear drift term and obtain $\fh \in \mcC([0,T]; L^1_1(\mR^d\times\mR^d))$ as in Definition \ref{def:weak solution}, we next prove the convergence in $L^1_1$.
    For any $R > 0$, we divide the domain into $\{|p|\le R\}$ and $\{|p| > R\}$.
    Since
    $$1+|p| \le \frac{(1+|p|)|p|}{R} \le \frac{2(1+|p|^2)}{R} \quad\text{for } |p|>R,$$
    we use \eqref{bound:unif fn 1} and \eqref{bound:unif fn 2} to obtain
    \begin{align*}
        \sup_{0\le t\le T}\|f_{\hbar,n}(t) - f_{\hbar,j}(t)\|_{L^1_1} &\le (1+R)\sup_{0\le t\le T}\|f_{\hbar,n}(t) - f_{\hbar,j}(t)\|_{L^1_{x,p}} \cr
        &\quad+ \frac{2}{R}\sup_{0\le t\le T}\inttr (1+|p|^2)\big(f_{\hbar,n}(t) + f_{\hbar,j}(t)\big)\,\rdx\rdp \cr
        &\le (1+R)\|f_{\hbar,0,n} - f_{\hbar,0,j}\|_{L^1_{x,p}} + \frac{4C_T}{R}.
    \end{align*}
    Taking $\limsup_{n,j\to\infty}$ first by \eqref{conv:init n} and then letting $R \to \infty$, we deduce that $\{f_{\hbar,n}\}$ is a Cauchy sequence in $\mcC([0,T]; L^1_1(\mR^d\times\mR^d))$.

    Thus, there exists a limit function $\fh \in \mcC([0,T]; L^1_1(\mR^d\times\mR^d))$ such that
    \begin{align}\label{conv:fn L11}
        \lim_{n\to\infty}\sup_{0\le t\le T}\|f_{\hbar,n}(t) - \fh(t)\|_{L^1_1} = 0.
    \end{align}
    From \eqref{conv:fn L11}, we immediately obtain
    $$\|\fh(t)\|_{L^1_{x,p}} = \lim_{n\to\infty}\|f_{\hbar,0,n}\|_{L^1_{x,p}} = \|f_{\hbar,0}\|_{L^1_{x,p}} \quad\text{for all } 0\le t\le T.$$
    Also, up to a subsequence, $f_{\hbar,n}$ converges to $\fh$ almost everywhere, so \eqref{bound:unif fn 1} and Fatou's lemma with \eqref{bound:unif fn 2} give
    $$0 \le \fh(x,p,t) \le \frac{1}{\hbar}\quad\text{and}\quad \fh \in L^\infty(0,T; L^1_2(\mR^d\times\mR^d)),$$
    satisfying \eqref{ineq:second_moment}.
    Since $\{\nabla_p f_{\hbar,n}\}$ is uniformly bounded in $L^2(0,T; L^2(\mR^d\times\mR^d))$ by \eqref{bound:unif fn 2}, we also obtain $\nabla_p \fh \in L^2(0,T; L^2(\mR^d\times\mR^d))$ by weak convergence.

    Next, $f_{\hbar,n}\in\Lambda_T$ satisfies \eqref{eq:int} for $f_{\hbar,0,n}$, and we have $\hbar|q||f_{\hbar,n}^2-\fh^2|\le 2|q||f_{\hbar,n}-\fh|$ from $0\le f_{\hbar,n},\fh\le 1/\hbar$.
    By Lemma \ref{Lem:F bound} with $k=l=1$ and $m=0$, we obtain for $0\le t\le T$:
    \begin{align*}
        &\lt\|f_{\hbar,n}(t) - \msF[f_{\hbar,0}](t) - \int_0^t \msF^q\big[\hbar q(\fh(s))^2\big](t-s)\,\rds\rt\|_{L^1_{x,p}} \cr
        &\quad \le C_1 e^{dt}\|f_{\hbar,0,n} - f_{\hbar,0}\|_{L^1_{x,p}} + C_2\int_0^t \frac{e^{d(t-s)}}{\sqrt{t-s}}\inttr \hbar|q|\big|f_{\hbar,n}(s)^2 - \fh(s)^2\big|\,\rdy\rdq\rds \cr
        &\quad \le C_1 e^{dT}\|f_{\hbar,0,n} - f_{\hbar,0}\|_{L^1_{x,p}} + 4C_2 e^{dT}\sqrt{T}\sup_{0\le s\le T}\|f_{\hbar,n}(s) - \fh(s)\|_{L^1_1}.
    \end{align*}
    Letting $n\to\infty$ with \eqref{conv:init n} and \eqref{conv:fn L11}, we see that $\fh$ satisfies the integral equation \eqref{eq:int}.

    Also, for $\phi \in \mcC_c^\infty([0,T) \times \mR^d \times \mR^d)$, we write \eqref{eq:weak form 1} for $f_{\hbar,n}$ as
    \begin{align}\label{eq:weak fn}
    \begin{split}
        &-\int_0^T \inttr f_{\hbar,n} (\partial_t \phi + p\cdot\nabla_x \phi + \Delta_p \phi - p \cdot \nabla_p \phi) \,\rdx\rdp\rdt \cr
        &\quad = \inttr f_{\hbar,0,n} \phi(x,p,0) \,\rdx\rdp + \int_0^T \inttr \hbar p f_{\hbar,n}^2 \cdot \nabla_p \phi \,\rdx\rdp\rdt.
    \end{split}
    \end{align}
    Note that since $\phi$ has compact support, $\frac{|\partial_t\phi + p\cdot\nabla_x \phi + \Delta_p \phi - p \cdot \nabla_p \phi|}{1+|p|}$, $|\phi(\cdot,\cdot,0)|$, and $|\nabla_p\phi|$ have an upper bound $C_\phi > 0$, and we have the bounds:
    \begin{align*}
        \lt|\int_0^T \inttr (f_{\hbar,n} - \fh)(\partial_t\phi + p\cdot\nabla_x \phi + \Delta_p \phi - p \cdot \nabla_p \phi)\,\rdx\rdp\rdt\rt| &\le C_\phi T\sup_{0\le t\le T}\|f_{\hbar,n}(t) - \fh(t)\|_{L^1_1}, \cr
        \lt|\inttr (f_{\hbar,0,n} - f_{\hbar,0})\phi(x,p,0)\,\rdx\rdp\rt| &\le C_\phi\|f_{\hbar,0,n} - f_{\hbar,0}\|_{L^1_{x,p}}, \cr
        \lt|\int_0^T \inttr \hbar p\big(f_{\hbar,n}^2 - \fh^2\big)\cdot\nabla_p\phi\,\rdx\rdp\rdt\rt| &\le 2C_\phi T\sup_{0\le t\le T}\|f_{\hbar,n}(t) - \fh(t)\|_{L^1_1}.
    \end{align*}
    Passing to the limit $n \to \infty$ in \eqref{eq:weak fn} with \eqref{conv:init n} and \eqref{conv:fn L11}, we see that $\fh$ satisfies \eqref{eq:weak form 1}.

    To end this proof, if $g_\hbar$ is another weak solution with the same initial data $f_{\hbar,0}$, then \eqref{ineq:contraction} gives
    $$\|\fh(t) - g_\hbar(t)\|_{L^1_{x,p}} \le \|f_{\hbar,0} - f_{\hbar,0}\|_{L^1_{x,p}} = 0.$$
\end{proof}

\subsection{Free energy inequality}

In this subsection, we prove the second part of Theorem \ref{theorem:exist} concerning the spatial second moment and the Fermi--Dirac free energy inequality \eqref{ineq:main_free_energy}.
For $0\le z\le 1/\hbar$, let
$$\eta_\hbar(z) := z\log z+\frac{1}{\hbar}(1-\hbar z)\log(1-\hbar z).$$
Since $\eta_\hbar'(z)=\log\frac{z}{1-\hbar z}$ for $0<z<1/\hbar$, we can write the drift and diffusion terms in \eqref{eq:quantum} as
$$\nabla_pg+p g(1-\hbar g)=g(1-\hbar g)\nabla_p\lt(\log\frac{g}{1-\hbar g}+\frac{|p|^2}{2}\rt)=g(1-\hbar g)\nabla_p\frac{\delta\mcH_\hbar}{\delta g},$$
which puts \eqref{eq:quantum} into the gradient-flow form in the momentum variable:
$$\partial_t\fh+p\cdot\nabla_x\fh=\nabla_p\cdot\lt[\fh(1-\hbar\fh)\nabla_p\frac{\delta\mcH_\hbar}{\delta\fh}\rt].$$

\begin{proof}[Proof of Theorem \ref{theorem:exist} (Free energy inequality)]
    We first derive the spatial second moment bound by using the cut-off functions $\zeta_N(x)$ and $\zeta_R(p)$ defined in Subsection \ref{sec:global}.
    Using $|x|^2\zeta_N(x)^2\zeta_R(p)$ as a test function in \eqref{eq:weak form 1} on $[0,t]$ and integrating by parts, we have for $t\in[0,T]$:
    \begin{align}\label{eq:moment_x_cutoff}
    \begin{split}
        &\inttr |x|^2\zeta_N(x)^2\zeta_R(p)\fh(x,p,t)\,\rdx\rdp - \inttr |x|^2\zeta_N(x)^2\zeta_R(p)f_{\hbar,0}(x,p)\,\rdx\rdp \cr
        &\quad = \int_0^t\inttr \fh\,p\cdot\nabla_x\big(|x|^2\zeta_N(x)^2\big)\zeta_R(p)\,\rdx\rdp\rds \cr
        &\qquad + \int_0^t\inttr |x|^2\zeta_N(x)^2\lt(\Delta_p\zeta_R(p)\,\fh - \nabla_p\zeta_R(p)\cdot p\fh(1-\hbar\fh)\rt)\,\rdx\rdp\rds.
    \end{split}
    \end{align}
    Since $0\le \zeta_N(x)\le 1$, $\supp(\nabla_x\zeta_N)\subset\{|x|\le 2N\}$, and $|\nabla_x\zeta_N(x)|\le C/N$, we have
    $$\big|\nabla_x\big(|x|^2\zeta_N(x)^2\big)\big| = \big|2x\zeta_N(x)^2 + 2|x|^2\zeta_N(x)\nabla_x\zeta_N(x)\big| \le C_0|x|\zeta_N(x),$$
    where $C_0>0$ is independent of $N$.
    By Young's inequality and $0\le \zeta_R(p)\le 1$, the first integral on the right-hand side of \eqref{eq:moment_x_cutoff} is bounded by
    $$\int_0^t\inttr \big|\fh\,p\cdot\nabla_x\big(|x|^2\zeta_N(x)^2\big)\zeta_R(p)\big|\,\rdx\rdp\rds \le \frac{C_0}{2}\int_0^t\inttr \big(|x|^2\zeta_N(x)^2 + |p|^2\big)\fh\,\rdx\rdp\rds.$$
    In the second integral on the right-hand side of \eqref{eq:moment_x_cutoff}, $|x|^2\zeta_N(x)^2\le 4N^2$, $0\le\fh(1-\hbar\fh)\le\fh$, $|\nabla_p\zeta_R(p)|\le C/R$, and $|\Delta_p\zeta_R(p)|\le C/R^2$, so for each fixed $N\ge 1$, \eqref{ineq:second_moment} and $\|\fh(s)\|_{L^1_{x,p}}=\|f_{\hbar,0}\|_{L^1_{x,p}}$ yield
    \begin{align*}
        &\lt|\int_0^t\inttr |x|^2\zeta_N(x)^2\lt(\Delta_p\zeta_R(p)\,\fh - \nabla_p\zeta_R(p)\cdot p\fh(1-\hbar\fh)\rt)\,\rdx\rdp\rds\rt| \cr
        &\quad \le 4N^2\int_0^t\inttr \lt(\frac{C}{R^2}\fh + \frac{C}{R}|p|\fh\rt)\,\rdx\rdp\rds \to 0 \quad\text{as } R\to\infty.
    \end{align*}
    When $R\to\infty$ with $N\ge 1$ fixed, \eqref{eq:moment_x_cutoff} becomes
    \begin{align*}
        \inttr |x|^2\zeta_N(x)^2\fh(x,p,t)\,\rdx\rdp &\le \inttr |x|^2 f_{\hbar,0}(x,p)\,\rdx\rdp + \frac{C_0}{2}\int_0^t\inttr |p|^2\fh(x,p,s)\,\rdx\rdp\rds \cr
        &\quad + \frac{C_0}{2}\int_0^t\inttr |x|^2\zeta_N(x)^2\fh(x,p,s)\,\rdx\rdp\rds.
    \end{align*}
    Note that for fixed $N\ge 1$, the quantity $\inttr |x|^2\zeta_N(x)^2\fh(t)\,\rdx\rdp \le 4N^2\|f_{\hbar,0}\|_{L^1_{x,p}}$ is finite and continuous on $[0,T]$.
    Applying Gr\"onwall's inequality with \eqref{ineq:second_moment} yields a bound independent of $N$, and passing to the limit $N\to\infty$ by the monotone convergence theorem gives for all $t \in [0,T]$,
    \begin{align}\label{ineq:moment_x}
        \inttr |x|^2\fh(x,p,t)\,\rdx\rdp \le C_T.
    \end{align}

    Next, let $\Phi_\hbar(g,p):=\log\frac{g}{1-\hbar g}+\frac{|p|^2}{2}$. 
    Since $\Phi_\hbar(g,p)$ is not bounded near $g=0$ and $g=1/\hbar$, we introduce a smooth function $\theta_K(r)$ for $K\ge 1$ as follows:
    $$\theta_K(r):=\int_0^r\zeta\lt(\frac{|s|}{K}\rt)\,\rds,\quad\text{where }\zeta\text{ is the cut-off function in Subsection \ref{sec:global}}.$$
    From the properties of $\zeta$, it follows directly that
    $$0\le\theta_K'(r)=\zeta(|r|/K)\le1,\quad \theta_K(r)=r\ \text{for}\ |r|\le K,\quad \text{and}\quad |\theta_K(r)|\le MK$$ 
    with $M:=\int_0^2\zeta(s)\,\rds$.
    Also, we have $\lim_{K\to\infty}\theta_K'(r)=1$ and $\lim_{K\to\infty}\theta_K(r)= r$ for all $r\in\mR$.

    Now, we define
    $$\Psi_{\hbar,K}(z,p):=\int_0^z\theta_K\lt(\Phi_\hbar(s,p)\rt)\,\rds.$$
    Since $|\theta_K(r)|\le|r|$, we obtain
    $$|\Psi_{\hbar,K}(z,p)|\le\int_0^z\lt|\log\frac{s}{1-\hbar s}+\frac{|p|^2}{2}\rt|\,\rds.$$
    For each fixed $K\ge 1$ and $|p|\le 2R$, the function $z\mapsto\theta_K(\Phi_\hbar(z,p))$ is constant near $z=0$ and $z=1/\hbar$, and thus smooth and bounded on $[0,1/\hbar]$.
    Since
    $$\fh\in L^\infty(0,T;L^\infty(\mR^d\times\mR^d))\quad\text{and}\quad \nabla_p\fh\in L^2(0,T;L^2(\mR^d\times\mR^d)),$$
    a standard regularization argument allows us to use $\theta_K(\Phi_\hbar(\fh,p))\zeta_N(x)\zeta_R(p)$ as a test function.
    Testing \eqref{eq:quantum} with $\theta_K(\Phi_\hbar(\fh,p))\zeta_N(x)\zeta_R(p)$ and integrating by parts, we use the identities
    $$\nabla_p\theta_K(\Phi_\hbar(\fh,p))=\theta_K' (\Phi_\hbar(\fh,p))\nabla_p\Phi_\hbar(\fh,p),\qquad \nabla_p\fh+p\fh(1-\hbar\fh)=\fh(1-\hbar\fh)\nabla_p\Phi_\hbar(\fh,p),$$
    and
    $$\lt(\nabla_p\fh+p\fh(1-\hbar\fh)\rt)\cdot\nabla_p\theta_K(\Phi_\hbar(\fh,p))=\theta_K'(\Phi_\hbar(\fh,p))\fh(1-\hbar\fh)|\nabla_p\Phi_\hbar(\fh,p)|^2\ge0$$
    to obtain
    \begin{align}\label{eq:truncated entropy}
    \begin{split}
        &\frac{\rd}{\rdt}\inttr\Psi_{\hbar,K}(\fh,p)\zeta_N(x)\zeta_R(p)\,\rdx\rdp\cr
        &\quad+\inttr\zeta_N(x)\zeta_R(p)\theta_K'(\Phi_\hbar(\fh,p))\fh(1-\hbar\fh)|\nabla_p\Phi_\hbar(\fh,p)|^2\,\rdx\rdp \cr
        &\ =E_{N,R,K}^{(1)}(t)+E_{N,R,K}^{(2)}(t),
    \end{split}
    \end{align}
    where
    \begin{align*}
        E_{N,R,K}^{(1)}(t)&:=\inttr p\cdot\nabla_x\zeta_N(x)\zeta_R(p)\Psi_{\hbar,K}(\fh,p)\,\rdx\rdp, \cr
        E_{N,R,K}^{(2)}(t)&:=-\inttr\zeta_N(x)\nabla_p\zeta_R(p)\cdot\lt(\nabla_p\fh+p\fh(1-\hbar\fh)\rt)\theta_K(\Phi_\hbar(\fh,p))\,\rdx\rdp.
    \end{align*}

    From $|\theta_K(\Phi_\hbar)| \le MK$, we have $|\Psi_{\hbar,K}(\fh,p)| \le MK\fh$, and $|\nabla_x\zeta_N(x)|\le C/N$ gives
    $$|E_{N,R,K}^{(1)}(t)|\le\frac{CMK}{N}\inttr |p|\fh\,\rdx\rdp,$$
    so by \eqref{ineq:second_moment},
    $$\int_0^t|E_{N,R,K}^{(1)}(s)|\,\rds\le\frac{C_{T,K}}{N}\to0\quad\text{as }N\to\infty.$$
    For $E_{N,R,K}^{(2)}(t)$, we use
    $$\theta_K(\Phi_\hbar(\fh,p))\nabla_p\fh = \nabla_p\big(\Psi_{\hbar,K}(\fh,p)\big) - p\int_0^{\fh}\theta_K'(\Phi_\hbar(s,p))\,\rds.$$
    By integration by parts in $p$ with $|\nabla_p\zeta_R(p)| \le C/R$, $|\Delta_p\zeta_R(p)| \le C/R^2$, and $0 \le \theta_K' \le 1$, we get
    \begin{align*}
        |E_{N,R,K}^{(2)}(t)| &\le \lt|\inttr \zeta_N(x)\Delta_p\zeta_R(p)\Psi_{\hbar,K}(\fh,p)\,\rdx\rdp \rt|\cr
        &\quad  + \lt|\inttr \zeta_N(x)\nabla_p\zeta_R(p)\cdot p\lt(\int_0^{\fh}\theta_K'(\Phi_\hbar(s,p))\,\rds - \fh(1-\hbar\fh)\theta_K(\Phi_\hbar(\fh,p))\rt)\,\rdx\rdp\rt| \cr
        &\le \frac{CMK}{R^2}\inttr \fh\,\rdx\rdp + \frac{C(MK+1)}{R}\inttr |p|\fh\,\rdx\rdp.
    \end{align*}
    Then \eqref{ineq:second_moment} gives
    $$\int_0^t|E_{N,R,K}^{(2)}(s)|\,\rds\le C_{T,K}(R^{-1}+R^{-2})\to0\quad \text{as}\ R\to\infty.$$
    Integrating \eqref{eq:truncated entropy} on $[0,t]$ and taking $N,R\to\infty$ by the dominated convergence theorem and Fatou's lemma, we have
    \begin{align}\label{ineq:energy inequality}
    \begin{split}   
        &\inttr\Psi_{\hbar,K}(\fh(t),p)\,\rdx\rdp+\int_0^t\inttr\theta_K'(\Phi_\hbar(\fh,p))\fh(1-\hbar\fh)|\nabla_p\Phi_\hbar(\fh,p)|^2\,\rdx\rdp\rds \cr
        &\quad\le\inttr\Psi_{\hbar,K}(f_{\hbar,0},p)\,\rdx\rdp.
    \end{split}
    \end{align}
    
    Now we check a lower bound of $\Psi_{\hbar,K}$ independent of $K$.
    For $r<0$, $0\le\zeta\le1$ implies $\theta_K(r)\ge r$.
    Since $\Phi_\hbar(s,p)\ge \log s$ for $0<s\le 1/\hbar$, we get for $0\le z\le1/\hbar$:
    $$\Psi_{\hbar,K}(z,p) \ge -\int^z_0\lt(\Phi_\hbar(s,p)\rt)_-\,\rds \ge -\int^z_0(-\log s)_+\,\rds \ge -z - (z\log z)_-.$$
    Let $W(x,p):=1+|x|^2+|p|^2$.
    Using $(z\log z)_- \le \frac{W}{2}z$ for $z \ge e^{-W/2}$ and $(z\log z)_- = z^{1/2}(-z^{1/2}\log z) \le e^{-W/4}$ for $0 < z < e^{-W/2}$, we obtain for $0\le z\le 1/\hbar$,
    $$\Psi_{\hbar,K}(z,p) \ge -\lt(1+\frac{W}{2}\rt)z - e^{-W/4}.$$
    Since $\lt(1+\frac{W}{2}\rt)\fh + e^{-W/4} \in L^1(\mR^d\times\mR^d)$ by \eqref{ineq:second_moment} and \eqref{ineq:moment_x}, we can apply Fatou's lemma as $K\to\infty$ to the left-hand side of \eqref{ineq:energy inequality}, while the right-hand side converges to $\mcH_\hbar(f_{\hbar,0})$ by the dominated convergence theorem.
    Therefore, letting $K\to\infty$ in \eqref{ineq:energy inequality} yields
    $$\mcH_\hbar(\fh(t))+\int_0^t\mcD_\hbar(\fh(s)) \,\rds\le\mcH_\hbar(f_{\hbar,0}).$$
\end{proof}

\section{Semi-classical limit and convergence rate}\label{sec:semi}

In this section, we prove the semi-classical limit and the quantitative convergence rate in Theorem \ref{theorem:semi} for the global weak solution $\fh$ of \eqref{eq:quantum} obtained in Theorem \ref{theorem:exist}.
Before proving Theorem \ref{theorem:semi}, we state the basic properties of the solution $f$ to the classical Fokker--Planck equation \eqref{eq:classic} by using the Green's function $\mcF$.

\begin{lemma}\label{Lem:classic_prop}
    Let $f_0 \in L^1_2(\mR^d\times\mR^d)$ be nonnegative. 
    Then for any $T>0$, the unique global weak solution $f$ of \eqref{eq:classic} on $[0,T]$ is given by $f(x,p,t) = \msF[f_0](x,p,t)$ in $\mcC([0,T]; L^1_j(\mR^d\times\mR^d))$ for all $j \in [0, 2)$, and satisfies for all $0\le t\le T$:
    \begin{enumerate}
        \item[(i)] $f(x,p,t) \ge 0$, $\|f(t)\|_{L^1_{x,p}} = \|f_0\|_{L^1_{x,p}}$, and if $g(t)=\msF[g_0](t)$ for another $g_0 \in L^1_2(\mR^d\times\mR^d)$, then
        \begin{align}\label{ineq:classic_contraction}
            \|f(t) - g(t)\|_{L^1_{x,p}} \le \|f_0 - g_0\|_{L^1_{x,p}}.
        \end{align}
        \item[(ii)] We have the second moment bound:
        \begin{align}\label{ineq:classic_moment}
            \inttr |p|^2 f(x,p,t)\,\rdx\rdp \le \inttr |p|^2 f_0(x,p)\,\rdx\rdp + 2d\|f_0\|_{L^1_{x,p}}t.
        \end{align}
    \end{enumerate}
\end{lemma}
\begin{proof}
    Since \eqref{eq:classic} is linear, by \eqref{prop} and standard linear theory \cite{VO90}, the weak solution of \eqref{eq:classic} is uniquely given by $f(t)=\msF[f_0](t)$.
    By \eqref{def:Green} and \eqref{prop:Green_mass}, we have $\mcF \ge 0$ and $\inttr \mcF(x,y,p,q,t)\,\rdx\rdp = 1$.
    Thus, Fubini's theorem gives $f(t)=\msF[f_0](t) \ge 0$, $\|f(t)\|_{L^1_{x,p}} = \|f_0\|_{L^1_{x,p}}$, and
    $$\|f(t) - g(t)\|_{L^1_{x,p}} \le \inttr \lt(\inttr \mcF(x,y,p,q,t)\,\rdx\rdp\rt)|f_0(y,q) - g_0(y,q)|\,\rdy\rdq = \|f_0 - g_0\|_{L^1_{x,p}}.$$
    Also, multiplying \eqref{eq:classic} by $|p|^2\zeta_N(x)\zeta_R(p)$ and integrating by parts as in Subsection \ref{sec:global} with $\hbar=0$, we obtain
    $$\frac{\rd}{\rdt}\inttr |p|^2 f(x,p,t)\,\rdx\rdp = 2d\inttr f(x,p,t)\,\rdx\rdp - 2\inttr |p|^2 f(x,p,t)\,\rdx\rdp \le 2d\|f_0\|_{L^1_{x,p}},$$
    which gives \eqref{ineq:classic_moment}.
    Finally, since $\msF[f_0] \in \mcC([0,T]; L^1_{x,p})$, combining the continuity in $L^1_{x,p}$ with the uniform second moment bound \eqref{ineq:classic_moment} on $\{|p|>R\}$ shows that
    $$f \in \mcC([0,T]; L^1_j(\mR^d\times\mR^d))\ \text{for all}\ j \in [0, 2).$$
\end{proof}

Note that by Theorem \ref{theorem:exist} and \eqref{ineq:second_moment}, the solution $\fh$ of \eqref{eq:quantum} also satisfies 
$$\fh \in \mcC([0,T]; L^1_{x,p}) \cap L^\infty(0,T; L^1_2(\mR^d\times\mR^d)),$$
so $\fh \in \mcC([0,T]; L^1_j(\mR^d\times\mR^d))$ for any $j \in [0, 2)$.

Now, we prove Theorem \ref{theorem:semi}.

\begin{proof}[Proof of Theorem \ref{theorem:semi}]
    Let us fix $T>0$, $j \in [0, 2)$, and $\epsilon>0$.
    For $f_0 \in L^1_2(\mR^d\times\mR^d)$, we choose $g_0 \in \mcC_c^\infty(\mR^d\times\mR^d)$ with $g_0 \ge 0$ such that
    \begin{align}\label{ineq:g0_approx}
        \|g_0 - f_0\|_{L^1_{x,p}} < \epsilon.
    \end{align}
    When $\hbar > 0$ is small enough so that $\hbar\|g_0\|_{L^\infty_{x,p}} \le 1$, we have $0 \le g_0(x,p) \le 1/\hbar$.
    Let $g_\hbar \in \Lambda_T$ and $g$ be the global solutions of \eqref{eq:quantum} and \eqref{eq:classic} for the same initial data $g_0$.

    From \eqref{ineq:second_moment} and \eqref{ineq:Lr_bound} for $r=3$ in Lemma \ref{Lem:reg_global}, $g_\hbar$ satisfies
    $$\sup_{0\le t\le T}\|g_\hbar(t)\|_{L^3_{x,p}} \le e^{\frac{2d}{3}T}\|g_0\|_{L^3_{x,p}} \quad\text{and}\quad \sup_{0\le t\le T}\|g_\hbar(t)\|_{L^1_2} \le \|g_0\|_{L^1_2} + 2d\|g_0\|_{L^1_{x,p}}T.$$
    By the Cauchy--Schwarz inequality, we have for $0\le t\le T$:
    \begin{align}\label{ineq:gh_nonlinear}
    \begin{split}
        \inttr |p|g_\hbar(x,p,t)^2\,\rdx\rdp &\le \lt(\inttr g_\hbar(x,p,t)^3\,\rdx\rdp\rt)^{\frac{1}{2}}\lt(\inttr |p|^2 g_\hbar(x,p,t)\,\rdx\rdp\rt)^{\frac{1}{2}} \cr
        &\le \|g_\hbar(t)\|_{L^3_{x,p}}^{\frac{3}{2}}\|g_\hbar(t)\|_{L^1_2}^{\frac{1}{2}} \le C_{T,g_0},
    \end{split}
    \end{align}
    where $C_{T,g_0} > 0$ does not depend on $\hbar$.
    
    In \eqref{eq:int}, we replace $\fh$ with $g_\hbar$, and by Lemma \ref{Lem:classic_prop} we have
    $$g_\hbar(t) = \msF[g_0](t) + \int_0^t \msF^q\big[\hbar q(g_\hbar(s))^2\big](t-s)\,\rds \quad\text{and}\quad g(t) = \msF[g_0](t).$$
    Using Lemma \ref{Lem:F bound} for $k=l=1$ and $m=0$ with \eqref{ineq:gh_nonlinear}, we obtain the bound for $0\le t\le T$ as
    \begin{align}\label{ineq:gh_g_diff}
    \begin{split}
        \|g_\hbar(t) - g(t)\|_{L^1_{x,p}} &\le C_2\int_0^t \frac{e^{d(t-s)}}{\sqrt{t-s}}\inttr \hbar|q|g_\hbar(y,q,s)^2\,\rdy\rdq\rds \cr
        &\le 2C_2 e^{dT}\sqrt{T}C_{T,g_0}\hbar.
    \end{split}
    \end{align}
    From \eqref{ineq:contraction}, \eqref{ineq:classic_contraction}, \eqref{ineq:g0_approx}, and \eqref{ineq:gh_g_diff}, we obtain
    \begin{align*}
        \sup_{0\le t\le T}\|\fh(t) - f(t)\|_{L^1_{x,p}} &\le \sup_{0\le t\le T}\|\fh(t) - g_\hbar(t)\|_{L^1_{x,p}} + \sup_{0\le t\le T}\|g_\hbar(t) - g(t)\|_{L^1_{x,p}} + \sup_{0\le t\le T}\|g(t) - f(t)\|_{L^1_{x,p}} \cr
        &\le \|f_{\hbar,0} - g_0\|_{L^1_{x,p}} + 2C_2 e^{dT}\sqrt{T}C_{T,g_0}\hbar + \|g_0 - f_0\|_{L^1_{x,p}} \cr
        &\le \|f_{\hbar,0} - f_0\|_{L^1_{x,p}} + 2\epsilon + 2C_2 e^{dT}\sqrt{T}C_{T,g_0}\hbar.
    \end{align*}
    Since $\epsilon>0$ can be chosen arbitrarily small, \eqref{assu:hbar initial} implies as $\hbar\to 0$ that
    \begin{align}\label{conv:fh_f_L1}
        \lim_{\hbar\to 0}\sup_{0\le t\le T}\|\fh(t) - f(t)\|_{L^1_{x,p}} = 0.
    \end{align}

    Now let $j\in[0,2)$ and $R \ge 1$.
    Since $1+|p|^j \le 1+R^j$ when $|p|\le R$ and $1+|p|^j \le \frac{2(1+|p|^2)}{R^{2-j}}$ when $|p|>R$, the second moment bounds \eqref{assu:moment}, \eqref{ineq:second_moment}, and \eqref{ineq:classic_moment} yield
    \begin{align*}
        \sup_{0\le t\le T}\|\fh(t) - f(t)\|_{L^1_j} &\le (1+R^j)\sup_{0\le t\le T}\|\fh(t) - f(t)\|_{L^1_{x,p}} \cr
        &\quad + \frac{2}{R^{2-j}}\sup_{0<\hbar\le 1}\sup_{0\le t\le T}\inttr (1+|p|^2)\big(\fh(t) + f(t)\big)\,\rdx\rdp \cr
        &\le (1+R^j)\sup_{0\le t\le T}\|\fh(t) - f(t)\|_{L^1_{x,p}} + \frac{2C_T}{R^{2-j}}.
    \end{align*}
    Taking $\limsup_{\hbar\to 0}$ first by \eqref{conv:fh_f_L1} and then letting $R \to \infty$ with $2-j > 0$, we get the convergence in $\mcC([0,T]; L^1_j(\mR^d\times\mR^d))$.

    Finally, when $f_0\in L^\infty(\mathbb R^d\times\mathbb R^d)$ and $\hbar\|f_0\|_{L^\infty_{x,p}}\le1$, we can take $g_0=f_0$ directly without $\epsilon$-approximation, so $g=f$ and $g_\hbar$ is the global weak solution of \eqref{eq:quantum} with initial datum $f_0$, satisfying the integral equation \eqref{eq:int}.
    Approximating $f_0$ by $f_{0,n} := f_0\mathbbm{1}_{\{|x|\le n, |p|\le n\}}$ as in the proof of Theorem \ref{theorem:exist} and letting $r \to \infty$ in \eqref{ineq:Lr_bound}, we obtain from \eqref{ineq:second_moment} and \eqref{ineq:Lr_bound} that
    $$\sup_{0\le t\le T}\|g_\hbar(t)\|_{L^\infty_{x,p}} \le e^{dT}\|f_0\|_{L^\infty_{x,p}} \quad\text{and}\quad \sup_{0\le t\le T}\|g_\hbar(t)\|_{L^1_2} \le \|f_0\|_{L^1_2} + 2d\|f_0\|_{L^1_{x,p}}T.$$
    Using $(1+|q|)|q| \le 2(1+|q|^2)$, we estimate the quadratic drift term on $[0,T]$ by
    \begin{align}\label{ineq:gh_L11_source}
    \begin{split}
        \big\|q(g_\hbar(s))^2\big\|_{L^1_1} &= \inttr (1+|q|)|q|g_\hbar(y,q,s)^2\,\rdy\rdq \cr
        &\le 2\|g_\hbar(s)\|_{L^\infty_{x,p}}\inttr (1+|q|^2)g_\hbar(y,q,s)\,\rdy\rdq \le C_{T,f_0},
    \end{split}
    \end{align}
    for a constant $C_{T,f_0} > 0$ depending only on $d, T, \|f_0\|_{L^1_2}$, and $\|f_0\|_{L^\infty_{x,p}}$.
    From Lemma \ref{Lem:F bound} with $k=l=1$ and $m=1$ and \eqref{ineq:gh_L11_source}, it follows that for all $t\in[0,T]$,
    \begin{align}\label{ineq:gh_f_L11_base}
        \|g_\hbar(t) - f(t)\|_{L^1_1} \le C_2\int_0^t \frac{e^{d(t-s)}}{\sqrt{t-s}}\hbar\big\|q(g_\hbar(s))^2\big\|_{L^1_1}\,\rds \le 2C_2 e^{dT}\sqrt{T}C_{T,f_0}\hbar.
    \end{align}
    For general $j \in [0, 2)$ and $R \ge 1$, we split the momentum integration into the regions $|p|\le R$ and $|p| > R$.
    Noting that $1+|p|^j \le (1+R^{(j-1)_+})(1+|p|)$ when $|p|\le R$ and $1+|p|^j \le \frac{2(1+|p|^2)}{R^{2-j}}$ when $|p|>R$, we combine \eqref{ineq:second_moment}, \eqref{ineq:classic_moment}, and \eqref{ineq:gh_f_L11_base} to deduce
    \begin{align*}
        \sup_{0\le t\le T}\|g_\hbar(t) - f(t)\|_{L^1_j} &\le (1+R^{(j-1)_+})\sup_{0\le t\le T}\|g_\hbar(t) - f(t)\|_{L^1_1} \cr
        &\quad + \frac{2}{R^{2-j}}\sup_{0\le t\le T}\inttr (1+|p|^2)\big(g_\hbar(t) + f(t)\big)\,\rdx\rdp \cr
        &\le 2C_2 e^{dT}\sqrt{T}C_{T,f_0}(1+R^{(j-1)_+})\hbar + \frac{4(\|f_0\|_{L^1_2} + 2d\|f_0\|_{L^1_{x,p}}T)}{R^{2-j}}.
    \end{align*}
    Setting $R = \max\{1, \hbar^{-1}\}$ yields the desired bound \eqref{ineq:quant_rate} for every $j \in [0, 2)$.
\end{proof}

\section*{Acknowledgements}
The author thanks Young-Pil Choi for his guidance and support, and Byung-Hoon Hwang for helpful discussions. 
This work received funding from NRF grant no. 2022R1A2C1002820 and no. RS-2024-00406821.

\bibliographystyle{abbrv}
\bibliography{FDFP_ref}

\end{document}